\documentclass[1p]{elsarticle}

\usepackage{amsmath}%
\usepackage{amsfonts}%
\usepackage{amssymb}%
\usepackage{graphicx}

\newtheorem{theorem1}{Theorem}
\newtheorem{lemma1}[theorem1]{Lemma}

\newtheorem{corollary1}[theorem1]{Corollary}
\newtheorem{definition1}[theorem1]{Definition}
\newtheorem{remark1}[theorem1]{Remark}

\def\Xint#1{\mathchoice
  {\XXint\displaystyle\textstyle{#1}}%
  {\XXint\textstyle\scriptstyle{#1}}%
  {\XXint\scriptstyle\scriptscriptstyle{#1}}%
  {\XXint\scriptscriptstyle\scriptscriptstyle{#1}}%
  \!\int}
\def\XXint#1#2#3{{\setbox0=\hbox{$#1{#2#3}{\int}$}
  \vcenter{\hbox{$#2#3$}}\kern-.5\wd0}}
\def\Mint{\Xint -}

\journal{Arxiv}

\begin{document}

\begin{frontmatter}

\title{A priori estimates for nonlinear fourth order Schr\"odinger type equations}

\author{G. M. Cardenas}
\ead{giovannomarcelorenato.cardenas@studenti.units.it}
\address[Dipartimento di Matematica]{Universit\'a Degli Studi di Trieste}

\begin{abstract}
We study the fourth order Schr\"odinger type differential inequality $-\Delta^2 u + \lambda V(x)u \geq a(x)u^q$ with $a,V\in L^1_{loc}(\mathbf{R}^N)$, both nonnegative, and $\lambda>0$.\\
We consider nonnegative solutions without making any assumptions about their asymptotic behavior, and choosing appropriate test functions, we find a priori estimates for some type of subharmonic solutions. Finally, applying these a priori estimates, we proof a Liouville type result.
\end{abstract}


\end{frontmatter}

\section{Introduction}
In this paper we will establish some Liouville type results for the fourth order Schr\"odinger type differential inequality
\begin{equation}\label{eq1.1} -\Delta^2 u + \lambda V(x)u \geq a(x)u^q, \qquad on \quad \mathbf{R}^N,\end{equation}
where $N>4$, $q>1$, $\lambda>0$ and $a, V \in L^1_{loc}(\mathbf{R}^N)$ are nonnegative functions such that there exists positive constants $C_1,C_2,\theta$, with $0\leq\theta<4$, for which we have that
\\
\begin{equation}\label{eq1.2}
V(x) \leq \frac{C_1}{|x|^4},\qquad
a(x) \geq \frac{C_2}{|x|^\theta},\qquad for\;large\;|x|.
\end{equation}
We assume that the potential $V(x)$ satisfies a Hardy type inequality. More precisely,
\begin{equation}\label{eq1.2.2}
\displaystyle{\lambda_H\int u^2V(x) \leq \int{(\Delta u)^2}\qquad \forall \; u\in W^{2,2}(\mathbf{R}^N)},
\end{equation} 
where $\lambda_H:=\left(\frac{N(N-4)}{4}\right)^2.$
\begin{remark1}\label{rem1.1} These assumptions hold in particular for $V(X)=\frac{1}{|x|^4}$ and $a(x)=\frac{1}{|x|^\theta}$, with $0\leq \theta<4$ and $N>4$.
\end{remark1}
Several Liouville results exist for second order coercive and non coercive elliptic differential equations. For instance, we have the following two well known results (More general results of this type, regarding also quasilinear operators, can be found in  \cite{mit5}, \cite{mit6}, \cite{far1},\cite{far2}, \cite{mit3}, \cite{serrin}).\\\\
\textbf{Theorem (Gidas, Spruck 1981).} 
Let $N>2$ and $1<q<\frac{N+2}{N-2}$. If $u\in C^2(\mathbf{R}^N)$ and $u>0$, then the equation
$$-\Delta u = u^q \quad on\;\mathbf{R}^N,$$
has no solution.\\\\
\textbf{Theorem (Brezis 1984).} Let $q>1$. If $u\in L^q_{loc}(\mathbf{R}^N)$ is a distributional solution of
\begin{equation}\label{eq1.5} \Delta u \geq |u|^{q-1} u \qquad on \; \mathbf{R}^N,\end{equation}
then $u\leq 0$ a.e. on $\mathbf{R}^N$. In particular if equality holds in (\ref{eq1.5}), then $u\equiv 0$ a.e. on $\mathbf{R}^N$.\\\\
An important thing to notice about the last theorem, is that no assumption regarding the asymptotic behavior of the solution is made. In fact, the result would be much simpler to proof otherwise.
In the present paper we will also make only local assumptions about the solutions. More precisely, we will consider solutions in the sense of the following definition.
\begin{definition1}\label{def1.1} Let $N>4$. A weak solution of (\ref{eq1.1}) is a function $u$ which satisfies\\\\
(i) $u\in W^{2,2}_{loc}(\mathbf{R}^N)\cap L^{\infty}_{loc}(\mathbf{R}^N)$\\
(ii) $\displaystyle{-\int\Delta u\Delta \phi + \int\lambda uV(x)\phi \geq \int u^q a(x)\phi \qquad \forall \phi\geq 0, \quad \phi\in C^{\infty}_0(\mathbf{R}^N)}$.
\end{definition1}
The fact that we don't have information about the asymptotic behavior of the solutions, makes finding a priori estimates an important issue. In order to do so, we will choose appropriate test functions by getting inspiration from \cite{mit3}.\\
This line of reasoning has already been applied for second order operators. In fact nonexistence results where established with a similar approach in \cite{mit1} for the following related differential inequality.
\begin{equation}\label{eq1.6} \Delta u + \lambda uV(x) \geq a(x)u^q\qquad on \quad \mathbf{R}^N.\end{equation}
Using this approach, a very precise meaning of the exponent of the nonlinearity was given, i.e. that the nonexistence results hold if and only if the exponent $q$ of the nonlinearity is less than some critical value. This makes us believe that developing further this technique, we could eventually arrive at a similar description of the nonlinearity exponent for the fourth order case.\\\\
We would like to mention that a similar problem to inequality (\ref{eq1.6}), was studied before in \cite{bly} using a different method. In that paper the authors reduced the problem to an ordinary differential inequality, and then applying some properties of second order ordinary differential equations, they proved nonexistence results.\\
One could think of adapting this ODE approach in order to deal with the fourth order case, but after reducing the problem to an ordinary differential inequality, it is not clear how to proceed, since to our knowledge, there are no analogous results for fourth order ordinary differential equations to the ones applied in that paper.\\\\
Together with conditions $(i)$ and $(ii)$ of Definition \ref{def1.1} of weak solution, we will also make assumptions about the sign of the solution and about the sign of its Laplacian.\\
In many problems involving quasilinear operators, no generality is lost when solutions are assumed to be positive, and this is due to the fact that there is a quasilinear version of Kato's inequality (for more details see \cite{mit2}). However, to our knowledge there is no analogous result for higher order operators. therefore more work needs to be done in order to treat the case when no assumptions are made about the sign of the solutions.\\\\
The plan of the paper is as follows. In Section 2 we give the proof of some preliminary results, regarding the properties of an auxiliary family of functions that will be used in the rest of the paper. Then, in Section 3, we obtain a priori estimates for the subharmonic solutions of (\ref{eq1.1}) and then proof their nonexistence.


\section{Preliminaries}
\subsection{Notation}
Let $R>0$. Throughout the paper we will assume that $A_R= \left\{x\in\mathbf{R}^N\; |\; R\leq |x|<2R\right\}$ and $B_R= \left\{x\in\mathbf{R}^N\; |\; |x|< R\right\}$. 
We denote by $C$ a generic nonnegative constant. The exact value of that constant may change from line to line and it will be specified whenever confusion may arise.\\
We will tacitly assume that $\int f(x)$ denotes the integral of $f$ over $\mathbf{R}^N$.

\subsection{Preliminary Results}
The following preliminary results describe some properties of an auxiliary family of functions that will be used in order to obtain a priori estimates of the possible solutions.\\
In what follows, we say that $\eta$ is a multi-index if $\eta = (\eta_1,\eta_2,...,\eta_N)$  where $\eta_1,...\eta_N$ are positive integers and we will use the notation $|\eta|:=\eta_1+...+\eta_N$. Also if $f\in C^\infty_0(\mathbf{R^N})$, we will denote $\frac{\partial^{|\eta|}f(x)}{(\partial x_1)^{\eta_1}(\partial x_i)^{\eta_2}...(x_N)^{\eta_N}}$ by $D^{\eta}f(x)$. Finally, if $f\in C^1(\mathbf{R})$, we will denote by $f'$ its derivative.\\

\begin{lemma1}\label{lem2.1}
Let $R>0$, $m>1$ and $\phi (x) := \left[\phi_0\left(\frac{|x|}{R}\right)\right]^m$ $\forall x=(x_1,x_2,...,x_N)\in\mathbf{R}^N$,  where
\begin{equation*}\phi_0(z) = \begin{cases}1 & 0 \leq |z| \leq 1\\
0 &|z|\geq2\end{cases},\end{equation*}
$\phi_0\in C^\infty_0(\mathbf{R})$ and $0<\phi_0(z)\leq 1$ for $0\leq z<2$.\\
Let $\eta$ be a multi-index and $0<|\eta|<m$, then\\
\begin{equation}\label{eq2.1} \;D^\eta\phi(x) = \sum_{s\in S} \left[\phi_0\left(\frac{|x|}{R}\right)\right]^{m-|\eta|}C_s\left(\frac{|x|}{R}\right)\frac{1}{|x|^{\beta_s}R^{\kappa_s}}\prod_{i=1}^N x_i^{\gamma_s^i},\qquad \forall x\in A_R\end{equation}
for some finite set of indexes $S$ which depend on $\eta$,  where $C_s\in C_0^\infty(\mathbf{R})$, $\gamma_s^i, \beta_s,\kappa_s\geq 0$ are integer that depend on the index $s\in S$, and
$$\quad \beta_s+\kappa_s-\sum_{i=1}^N\gamma_s^i = |\eta| \qquad \forall s\in S.$$
Also, if $d>0$, we obtain\\
\begin{equation} \label{eq2.2}\left|D^\eta\phi\right(x)|^d \leq C\frac{1}{R^{|\eta|d}} \left[\phi_0\left(\frac{|x|}{R}\right)\right]^{d(m-|\eta|)} \qquad \forall x\in A_R.\end{equation}
For some positive constant $C$ which does not depend on $R$.\\
\textbf{Proof.} We have that
$$\frac{\partial}{\partial x_i}\phi(x)=m\left[\phi_0\left(\frac{|x|}{R}\right)\right]^{m-1}\phi_0'\left(\frac{|x|}{R}\right)\frac{x_i}{|x|R}, \qquad 1\leq i\leq N,$$
and this shows that the statement is true when $|\eta|=1$.\\ Assuming that (\ref{eq2.1}) is true for some fixed $\eta=(\eta_1,\eta_2,...,\eta_N)$, with $m>|\eta|+1$, and defining $\hat{\eta}:=(\eta_1,\eta_2,...,\eta_j+1,...,\eta_N)$, we have that,
\begin{equation*}
\begin{array}{l}
\displaystyle{D^{\hat\eta}\phi(x)=\frac{\partial}{\partial x_j}(D^\eta\phi(x)) =}\\
\displaystyle{\sum_{s\in S}(m-|\eta|)\left[\phi_0\left(\frac{|x|}{R}\right)\right]^{m-|\eta|-1}\phi_0'\left(\frac{|x|}{R}\right)\frac{x_j}{|x|R}C_s\left(\frac{|x|}{R}\right)\frac{1}{|x|^{\beta_s}R^{\kappa_s}}\prod_{i=1}^N x_i^{\gamma_s^i}}\\
\displaystyle{+ \sum_{s\in S}\left[\phi_0\left(\frac{|x|}{R}\right)\right]^{m-|\eta|}C_s'\left(\frac{|x|}{R}\right)\frac{x_j}{|x|R}\frac{1}{|x|^{\beta_s}R^{\kappa_s}}\prod_{i=1}^N x_i^{\gamma_s^i}}\\
\displaystyle{+\sum_{s\in S}\left[\phi_0\left(\frac{|x|}{R}\right)\right]^{m-|\eta|}C_s\left(\frac{|x|}{R}\right)(-\beta_s)\frac{1}{|x|^{\beta_s+1}R^{\kappa_s}}\frac{x_j}{|x|}\prod_{i=1}^N x_i^{\gamma_s^i}}\\
\displaystyle{+ \sum_{s\in S} \left[\phi_0\left(\frac{|x|}{R}\right)\right]^{m-|\eta|}C_s\left(\frac{|x|}{R}\right)\frac{1}{|x|^{\beta_s}R^{\kappa_s}}\frac{\partial}{\partial x_j}\left(x_j^{\gamma_s^j}\right)\prod_{i\ne j} x_i^{\gamma_s^i},}\\
\end{array}
\end{equation*}
which means that,
\begin{equation*}
\begin{array}{l}
\displaystyle{D^{\hat\eta}\phi(x)= \left[\phi_0\left(\frac{|x|}{R}\right)\right]^{m-|\eta|-1}\prod_{i\ne j}^N x_i^{\gamma_s^i} \sum_{s\in S} C_s^1\left(\frac{|x|}{R}\right)\frac{x_j^{\gamma_S^j + 1}}{|x|^{\beta_s+1}R^{\kappa_s+1}}}\\
\displaystyle{+ C_s^2\left(\frac{|x|}{R}\right)\frac{x_j^{\gamma_s^j +1}}{|x|^{\beta_s + 2}R^{\kappa_s}} + C_s^3\left(\frac{|x|}{R}\right)\frac{x_j^{\gamma_s^j -1}}{|x|^{\beta_s}R^\kappa_s}
},\\
\end{array}
\end{equation*}
where $C_s^1,C_s^2,C_s^3\in C_0^\infty(\mathbf{R})$ depend on $\phi_0,\phi_0' , C_s, C_s'$ $\forall s\in S$. Thus the first part of the lemma follows from induction. Inequality (\ref{eq2.2}) easily follows from (\ref{eq2.1}).$\;\Diamond$\\
\end{lemma1}

\begin{corollary1}\label{cor2.1}
Let $\alpha>0$, $R_0>0$, $1<q<\frac{N - \theta + \alpha(4 - \theta)}{N-4}$, $\chi = \frac{\alpha + q}{\alpha + 1}$ and $\chi'=\frac{\alpha+q}{q-1}$ its conjugate exponent. Let $\phi$ be a family of functions depending on a positive real parameter $R$, obtained from the $\phi$ defined in the previous Lemma by leaving $R$ as a parameter and fixing a large enough exponent $m$ (Large enough such that the inequalities appearing in (\ref{eq2.3}) and (\ref{eq2.4}) hold). Then, we get the following.\\
(i) $\displaystyle{\lim_{R\to\infty}\int_{A_R}{|x|^{\frac{\theta\chi'}{\chi}}\frac{|\nabla \phi|^{4\chi'}}{\phi^{3\chi' + \frac{\chi'}{\chi}}} = 0}},$\\
(ii)$\displaystyle{\lim_{R\to\infty}\int_{A_R}{|x|^{\frac{\theta\chi'}{\chi}}\frac{(\Delta \phi)^{2\chi'}}{\phi^{\chi'+ \frac{\chi'}{\chi}}} =0}},$\\
(iii)$\displaystyle{\lim_{R\to\infty}\int^R_{R_0}{|x|^{\frac{\theta\chi'}{\chi}}|x|^{-4\chi'}\phi<\infty}}.$\\
\textbf{Proof.} Applying (\ref{eq2.2}) of Lemma \ref{lem2.1} to (i) and (ii), we obtain that for every $R>0$, the following holds.\\
\begin{equation}\label{eq2.3}
\displaystyle{\int_{A_R}{|x|^\frac{\theta\chi'}{\chi}\frac{|\nabla \phi|^{4\chi'}}{\phi^{3\chi'+ \frac{\chi'}{\chi}}} \leq C\int_{A_R}{\frac{|x|^\frac{\theta\chi'}{\chi}}{R^{4\chi'}\phi^{3\chi'+ \frac{\chi'}{\chi}}}  \left[\phi_0\left(\frac{|x|}{R}\right)\right]^{4(m-1)\chi'}} \leq CR^{\frac{\theta\chi'}{\chi}-4\chi' + N}}}\qquad\end{equation}
and also,
\begin{equation}\label{eq2.4} \displaystyle{\int_{A_R}{|x|^{\frac{\theta\chi'}{\chi}}\frac{(\Delta \phi)^{2\chi'}}{\phi^{\chi'+\frac{\chi'}{\chi}}} \leq C\int_{A_R}{\frac{|x|^\frac{\theta\chi'}{\chi}}{R^{4\chi'}\phi^{\chi' + \frac{\chi'}{\chi}}}  \left[\phi_0\left(\frac{|x|}{R}\right)\right]^{2(m-2)\chi'}} \leq CR^{\frac{\theta\chi'}{\chi}-4\chi' + N}}}\end{equation}
Moreover, for (iii) we have the following estimate.
$$\displaystyle{\int_{R_0}^R{|x|^{\frac{\theta\chi'}{\chi}}|x|^{-4\chi'}} =}$$
\begin{equation}\label{eq2.5} \displaystyle{C\int_{R_0}^R t^{\frac{\theta\chi'}{\chi} - 4\chi'}t^{N-1}dt \leq \frac{C}{\frac{\theta\chi'}{\chi}-4\chi' + N}\left[R^{\frac{\theta\chi'}{\chi}-4\chi' + N} - {R_0}^{\frac{\theta\chi'}{\chi}-4\chi' + N}\right] },\end{equation}
\\
where the constants C do not depend on $R$.\\ 
We also know that
$$\frac{\theta\chi'}{\chi}-4\chi'+N=\frac{\theta(1+\alpha)}{q-1} - \frac{4(q+\alpha)}{q-1} + N,$$ 
and since
$$q<\frac{N - \theta + \alpha(4 - \theta)}{N-4},$$
the claim holds because this last two inequalities imply that the exponent of $R$ in the r.h.s of (\ref{eq2.3}),(\ref{eq2.4}) and (\ref{eq2.5}) is negative. $\;\Diamond$
\end{corollary1}

\begin{remark1}\label{rem2.1} 
From now on, we will denote by $\phi$ a family of functions depending on a positive real parameter $R$, obtained from the $\phi$ defined in Lemma \ref{lem2.1} by fixing an appropiate exponent $m$ and leaving $R$ as a parameter. We will not write this dependency on $R$ explicitly, and every time we use a family $\phi$ in a proof where constants $\alpha$ and $q$ are fixed (as in the previous Corollary), we will tacitly assume that it is a family defined with exponent $m$ large enough, such that the previous Corollary holds for those values of $q$ and $\alpha$. \end{remark1}
It is easy to show by a regularization argument that if  $f\in W^{2,2}_{loc}(\mathbf{R}^N)$ is a nonnegative function, then the family $f\phi$ can be used as a test function in (ii) of Definition \ref{def1.1}. More precisely, we have the following lemma.

\begin{lemma1}\label{lem2.2} Let $u$ be a weak solution of (\ref{eq1.1}) and let $v\in W_{loc}^{2,2}(\mathbf{R}^N)$ be a nonnegative function. Then,
\begin{equation*}\label{eq2.5.3}-\displaystyle{\int\Delta u\Delta(v\psi) + \lambda\int uvV(x)\psi\geq \int u^qva(x)\psi}, \qquad  \forall \psi\geq 0, \; \psi\in C^{\infty}_0(\mathbf{R}^N).\end{equation*}
\textbf{Proof.} Let $(\rho_n)_{n\geq0}$ be a regularizing sequence and define $v_n:= (v\psi) * \rho_n\in C_0^{\infty}({\mathbf{R}^N})$. Also, let $\Omega$ be an open ball in $\mathbf{R}^N$, such that $supp(v_n) \bigcup supp(\psi) \subset \Omega$ $\forall n>0$. Since $v\psi|_{\Omega}\in W_0^{2,2}(\Omega)$, there exists a subsequence of $(v_n)_{n>0}$, which we again denote by $(v_n)_{n>0}$, such that
\\
\begin{equation}\label{eq2.5.4}
\begin{array}{ll}
	v_n\rightarrow v\psi & in \; W^{2,2}(\Omega)\\
 	v_n\rightarrow v\psi & a.e. \; in\; \Omega.
\end{array}
\end{equation}\\
From Hardy's inequality (\ref{eq1.2.2}) it follows that
$$\displaystyle{\lambda_H\int_{\Omega} |v_n-v_m|^2V(x)\leq \int_{\Omega} \left[\Delta(v_n-v_m)\right]^2}.$$
Therefore, $(v_n\left(V(x)\right)^\frac{1}{2})_{n\geq0}$ is a Cauchy sequence in $L^2(\Omega)$ and this means that we can obtain a further subsequence of $(v_n)_{n>0}$, which we continue to denote the same way, such that there exists $w\in L^2(\Omega)$ for which,
\begin{equation}\label{eq2.5.5}\;v_n\left(V(x)\right)^\frac{1}{2}\leq w\qquad  a.e.\; in\; \Omega \quad \forall\;n>0.\end{equation}
Moreover, applying Cauchy Schwartz inequality we also have that,
$$\int_{\Omega}|v_n-v|V(x)u\leq\left[\int_{\Omega}|v_n-v|^2V(x)\right]^\frac{1}{2}\left[\int_{\Omega}u^2V(x)\right]^\frac{1}{2},$$
and since (\ref{eq2.5.4}), (\ref{eq2.5.5}) and Lebesgue's Dominated convergence Theorem imply that,
$$v_n\left(V(x)\right)^\frac{1}{2}\rightarrow v\left(V(x)\right)^\frac{1}{2}\psi \qquad in \quad L^2(\mathbf{R}^N),$$
we can conclude that
\begin{equation}\label{eq2.5.8}\displaystyle{\int_{\Omega} v_nuV(x)\rightarrow\int_{\Omega} vuV(x)\psi},\end{equation}
using the facts that $u\in L^\infty_{loc}(\Omega)$ and $V\in L^1_{loc}(\mathbf{R}^N).$\\
Afterwards, we apply Cauchy-Schwarz inequality to obtain,
$$\displaystyle{\left|\int_{\Omega}\Delta u\Delta v_n - \int_{\Omega} \Delta u \Delta(v\psi)\right|\leq \left[\int_{\Omega}(\Delta u)^2\right]^{\frac{1}{2}}\left[\int_{\Omega}{|\Delta(v_n - v\psi)|^2}\right]^{\frac{1}{2}}},$$
and since $u\in W_{loc}^{2,2}(\mathbf{R}^N)$, it follows from (\ref{eq2.5.4}) that
\begin{equation}\label{eq2.5.7}\displaystyle{\int_{\Omega}\Delta u\Delta(v_n) \rightarrow \int_{\Omega}\Delta u \Delta(v\psi)}.\end{equation}
Finally, since $v_n\in C_0^{\infty}(\mathbf{R}^N)$, we get from (ii) of Definition \ref{def1.1} that
\begin{equation}\label{eq2.5.6}
\;-\displaystyle{\int\Delta u \Delta(v_n) + \lambda\int v_nuV(x)\geq \int v_nu^qa(x)\qquad \forall\;n>0}.
\end{equation}
Therefore, applying (\ref{eq2.5.8}) and (\ref{eq2.5.7}) we obtain that,\\
$$\displaystyle{\liminf\int v_nu^qa(x)\leq \lambda\int uvV(x)\psi - \int\Delta u\Delta(v\psi)},$$\\
and the proof follows after applying Fatou's Lemma to the last inequality. $\;\Diamond$
\end{lemma1}

\section{Subharmonic solutions}

\begin{definition1}\label{def3.1} We say that $u$ is a weak subharmonic solution of (\ref{eq1.1}) if $u$ is a weak solution in the sense of Definition $\ref{def1.1}$ and  $u\geq 0$, $\Delta u\geq 0$ a.e. in $\mathbf{R}^N.$
\end{definition1}
The following lemma contains the fundamental a priori estimates for the subharmonic solutions of (\ref{eq1.1}) .

\begin{lemma1}\label{lem3.1}
Let $u$ be a weak subharmonic solution of (\ref{eq1.1}) in the sense of Definition \ref{def3.1}. Let $\alpha>2$ and $1<q<\frac{N+\alpha(4-\theta)-\theta}{N-4}$, then\\
(i) $\;u^{q+\alpha}a(x)$, $u^{\alpha+1}V(x)$ $\in L^1(\mathbf{R}^N),$\\
(ii)$\;\Delta u\Delta(u^\alpha)\;$ $\in L^1(\mathbf{R}^N)$ and  $\;\displaystyle{\int u^{\alpha+q}a(x)+\int\Delta (u^{\alpha})\Delta u\leq \lambda\int{u^{\alpha+1}V(x)}}.$\\
\textbf{Proof.} (i) It is easy to see that $u^\alpha \in W_{loc}^{2,2}(\mathbf{R}^N)$ (see Appendix 5.1 for details). Therefore, applying Lemma \ref{lem2.2} we obtain,
$$\displaystyle{\int u^{\alpha+q}a(x)\phi \leq  -\int\Delta u\Delta(u^{\alpha}\phi)+\lambda\int u^{\alpha+1}V(x)\phi}.$$    
Then, because of the local properties of $u$ assumed in (i) of Definition \ref{def1.1}, it is possible to expand the Laplacian appearing in the first integral of the r.h.s. (see Appendix 5.1 for details) in order to obtain,
\begin{equation}\label{eq2.6}
\begin{array}{l}
\displaystyle{\int u^{\alpha+q}a(x)\phi+ \int\alpha u^{\alpha-1}(\Delta u)^2\phi + \int\alpha(\alpha-1)u^{\alpha-2}\Delta u|\nabla u|^2\phi} \\
  \displaystyle{ \leq -2\alpha\int{u^{\alpha-1}\Delta u\nabla u\nabla\phi} - \int u^{\alpha}\Delta u\Delta \phi + \lambda\int u^{\alpha+1}V(x)\phi}.
\end{array}
\end{equation}\\
We now look for upper bounds of the integrals in the r.h.s of the last inequality.\\
First, applying Young's inequality with exponent $\frac{1}{2}$ and parameter $\epsilon_1>0$ to the first integral appearing in the r.h.s. of (\ref{eq2.6}), we get that,
\begin{equation*}\displaystyle{2\alpha\left|\int u^{\alpha-1}\Delta u\nabla u\nabla\phi\right|\leq \frac{\epsilon_1^2}{2}\int\alpha(\alpha-1) u^{\alpha-2}\Delta u|\nabla u|^2\phi + \frac{4\alpha^2}{2\epsilon_1^2\alpha(\alpha-1)}\int u^{\alpha}\Delta u\frac{|\nabla\phi|^2}{\phi}},\end{equation*}
from which we obtain, applying Young's inequality with exponent $\frac{1}{2}$ and parameter $\epsilon_2>0$, that,
\begin{equation}\label{eq2.7} 
\begin{array}{l}
 2\displaystyle{\alpha\left|\int u^{\alpha-1}\Delta u\nabla u\nabla\phi\right|\leq}\\ \displaystyle{\frac{\epsilon_1^2}{2}\int\alpha(\alpha-1) u^{\alpha-2}|\nabla u|^2\Delta u\phi + \frac{\epsilon_2^2}{2}\int\alpha u^{\alpha-1}(\Delta u)^2\phi + k_1(\alpha)\int u^{\alpha+1}\frac{|\nabla \phi|^4}{\phi^3}},
\end{array}
\end{equation}
with $k_1=[\frac{4\alpha^2}{2\epsilon_1^2\alpha(\alpha-1)}]^2\frac{1}{2\alpha \epsilon_2^2}$.\\
Then, we apply Young's inequality with exponent $\frac{1}{2}$ and parameter $\epsilon_3>0$ to the second integral appearing in the r.h.s. of (\ref{eq2.6}),  and we get,
\begin{equation}\label{eq2.8}
\displaystyle{\left|\int u^{\alpha}\Delta u\Delta\phi\right|\leq \frac{\epsilon_3^2}{2}\int\alpha u^{\alpha-1}(\Delta u)^2\phi + \frac{1}{2\epsilon_3^2\alpha}\int u^{\alpha+1}\frac{(\Delta\phi)^2}{\phi}}.
\end{equation}
Now we apply the upper bounds obtained in (\ref{eq2.7}) and (\ref{eq2.8}) to inequality (\ref{eq2.6}), and this leads as to, 
\begin{equation}\label{eq2.9}
\begin{array}{l}
\displaystyle{\int u^{\alpha+q}a(x)\phi + k_3(\alpha)\int\alpha(\alpha-1)u^{\alpha-2}\Delta u|\nabla u|^2\phi + k_2(\alpha)\int\alpha u^{\alpha-1}(\Delta u)^2\phi \leq} \\ 
\displaystyle{k_1(\alpha)\int u^{\alpha+1}\frac{|\nabla\phi|^4}{\phi^3} +  \frac{1}{2\epsilon_3^2\alpha}\int u^{\alpha+1}\frac{(\Delta\phi)^2}{\phi}  +\lambda\int u^{\alpha+1}V(x)\phi},
\end{array}
\end{equation}
where $k_2(\alpha)=[1 - \frac{\epsilon_2^2}{2} - \frac{\epsilon_3^2}{2}]$, $k_3(\alpha)=1-\frac{\epsilon_1^2}{2}.$\\\\
Let now $\chi=\frac{q+\alpha}{\alpha+1}$ and $\chi'=\frac{q+\alpha}{q-1}$ be its conjugate exponent and let also $R_0>0$ be such that (\ref{eq1.2}) holds for $|x|\geq R_0$. Applying H\"older's inequality with exponent $\chi$ together with (\ref{eq1.2}), we obtain the following estimates for the integrals appearing in the r.h.s. of (\ref{eq2.9}). 
\begin{equation}\label{eq2.10}
\begin{array}{ll}
 (i') & \displaystyle{\left|k_1(\alpha)\int u^{\alpha+1}\frac{|\nabla\phi|^4}{\phi^3}\right| \leq k_1(\alpha)\left[\int\frac{u^{q+\alpha}}{|x|^{\theta}}\phi\right]^{\frac{1}{\chi}}\left[\int|x|^{\frac{\theta\chi'}{\chi}}\frac{|\nabla\phi|^{4\chi'}}{\phi^{3\chi'+\frac{\chi'}{\chi}}}\right]^{\frac{1}{\chi'}}},\\

(ii') & \displaystyle{\left|\frac{1}{2\epsilon_3^2\alpha}\int u^{\alpha+1}\frac{(\Delta \phi)^2}{\phi}\right|\leq\frac{1}{2\epsilon_3^2\alpha}\left[\int\frac{u^{q+\alpha}}{|x|^{\theta}}\phi\right]^{\frac{1}{\chi}}\left[\int|x|^{\frac{\theta\chi'}{\chi}}\frac{|\Delta\phi|^{2\chi'}}{\phi^{\chi'+\frac{\chi'}{\chi}}}\right]^{\frac{1}{\chi'}}},\\

(iii') & \displaystyle{\lambda\left|\int u^{\alpha+1}V(x)\phi\right|\leq \lambda\int_{|x|<R_0}{u^{\alpha+1}V(x)\phi} + C_1\lambda \int_{|x|\geq R_0}\frac{u^{1+\alpha}}{|x|^4}\phi \leq} \\
 & \displaystyle{\lambda\int_{|x|<R_0}{u^{\alpha+1}V(x)\phi} + C_1\lambda\left[\int_{|x|\geq R_0}\frac{u^{q+\alpha}}{|x|^{\theta}}\phi\right]^{\frac{1}{\chi}}\left[\int_{|x|\geq R_0}|x|^{\frac{\theta\chi'}{\chi}}|x|^{-4\chi'}\phi\right]^{\frac{1}{\chi'}}}.
\newline
\end{array}
\end{equation}
\\
Furthermore, from an application of Young's inequality with exponent $\chi$ and parameters $\delta_1>0,\delta_2>0$ and $\delta_3>0$ to (i'),(ii') and (iii') respectively, followed by an application of (\ref{eq1.2}), we get,\\
\begin{equation}\label{eq2.11}
\begin{array}{ll}

(i'') & \displaystyle{\left|k_1(\alpha)\int u^{\alpha+1}\frac{|\nabla\phi|^4}{\phi^3}\right| \leq \frac{\delta_1^{\chi}}{C_2\chi}\int u^{\alpha+q}a(x)\phi + \frac{[k_1(\alpha)]^{\chi'}}{\chi'\delta_1^{\chi'}}\left[\int|x|^{\frac{\theta\chi'}{\chi}}\frac{|\nabla\phi|^{4\chi'}}{\phi^{3\chi'+\frac{\chi'}{\chi}}}\right]},\\

(ii'') & \displaystyle{\left|\frac{1}{2\epsilon_3^2\alpha}\int u^{\alpha+1}\frac{(\Delta \phi)^2}{\phi}\right| \leq \frac{\delta_2^{\chi}}{C_2\chi}\int u^{\alpha+q}a(x)\phi + \left[\frac{1}{2\epsilon_3^2\alpha}\right]^{\chi'}\frac{1}{\chi'\delta_2^{\chi'}}\left[\int|x|^{\frac{\theta\chi'}{\chi}}\frac{|\Delta\phi|^{2\chi'}}{\phi^{\chi'+\frac{\chi'}{\chi}}}\right]},\\

(iii'') & \displaystyle{\lambda\left|\int u^{\alpha+1}V(x)\phi\right| \leq \lambda\int_{|x|<R_0}{u^{\alpha+1}V(x)\phi} + \frac{\delta_3^{\chi}}{\chi} \int_{|x|\geq R_0}\frac{u^{q+\alpha}}{|x|^\theta}\phi+  \frac{\lambda^{\chi'}C_1^{\chi'}}{\chi'\delta_3^{\chi'}}\int_{|x|\geq R_0}|x|^{\frac{\theta\chi'}{\chi}-4\chi'}\phi\leq}\\
& \displaystyle{\lambda\int_{|x|<R_0}{u^{\alpha+1}V(x)\phi}  + \frac{\delta_3^{\chi}}{C_2\chi} \int_{|x|\geq R_0}u^{q+\alpha}a(x)\phi+  \frac{\lambda^{\chi'}C_1^{\chi'}}{\chi'\delta_3^{\chi'}}\int_{|x|\geq R_0}|x|^{\frac{\theta\chi'}{\chi}-4\chi'}\phi}.
\end{array}
\end{equation}
\\
Applying (i''),(ii''),(iii'') to (\ref{eq2.9}) we obtain,\\
\\
\begin{equation}\label{eq2.12}
\begin{array}{l}
\displaystyle{
c_4(\alpha)\int u^{q+\alpha}a(x)\phi+ k_3(\alpha)\int\alpha(\alpha-1)u^{\alpha-2}\Delta u|\nabla u|^2\phi+ k_2(\alpha)\int\alpha u^{\alpha-1}(\Delta u)^2\phi}\\
\displaystyle{\leq \lambda\int_{|x|< R_0}u^{1+\alpha}V(x)\phi+ c_1(\alpha)\int_{A_R}|x|^{\frac{\theta\chi'}{\chi}}\frac{|\nabla\phi|^{4\chi'}}{\phi^{3\chi'+\frac{\chi'}{\chi}}} + c_2(\alpha)\int_{A_R}|x|^{\frac{\theta\chi'}{\chi}}\frac{|\Delta\phi|^{2\chi'}}{\phi^{\chi'+\frac{\chi'}{\chi}}}} \\
\displaystyle{+ c_3(\alpha)\int_{|x|\geq R_0}|x|^{\frac{\theta\chi'}{\chi}-4\chi'}\phi},
\end{array}
\end{equation}
\\
where,
$$
\begin{array}{ll}
c_1(\alpha)= \frac{[k_1(\alpha)]^{\chi'}}{\chi'\delta_1^{\chi'}}, & c_2(\alpha)= [\frac{1}{2\epsilon_3^2\alpha}]^{\chi'}\frac{1}{\chi'\delta_2^{\chi'}},\\\\
c_3(\alpha)= \frac{\lambda^{\chi'}C_1^{\chi'}}{\chi'\delta_3^{\chi'}} & c_4(\alpha)=(1-\frac{\delta_1^\chi}{C_2\chi}- \frac{\delta_2^\chi}{C_2\chi} - \frac{\delta_3^\chi}{C_2\chi}).
\end{array}
$$
Everything that has been said so far, holds for any choice of the positive parameters $\epsilon_1$, $\epsilon_2$, $\epsilon_3$, $\delta_1$, $\delta_2$, $\delta_3$, therefore, we may assume that they are small enough such that $c_4(\alpha),k_2(\alpha)$ and $k_3(\alpha)$ are positive, or equivalently, such that all the terms in the l.h.s. of  (\ref{eq2.12})  are nonnegative. Then, letting $R\rightarrow\infty$ we get from Corollary \ref{cor2.1} together with the Monotone Convergence Theorem, that 
$u^{q+\alpha}a(x) \in L^1(\mathbf{R}^N)$. Moreover, letting $R/rightarrow/infty$ in (iii') and applying Corollary \ref{cor2.1}, we get that $u^{\alpha +1}V(x)\in L^1(\mathbf{R}^N)$ and this proves (i).\\\\
In order to prove (ii) we assume that the positive parameters $\epsilon_1$, $\epsilon_2$ and $\epsilon_3$ are such that $1\geq k_3(\alpha)\geq k_2(\alpha)>0$ and are small enough such that everything said until now remains true. Then, applying (i') and (ii') from (\ref{eq2.10}) to (\ref{eq2.9}) we obtain,
\begin{equation}\label{eq2.12.1}
\begin{array}{l}
\;\displaystyle{\int u^{q+\alpha}a(x)\phi + k_2(\alpha)\int{\Delta u\Delta(u^\alpha)\phi} \leq}\\
\displaystyle{\int u^{q+\alpha}a(x)\phi + k_3(\alpha)\int\alpha(\alpha-1)u^{\alpha-2}\Delta u|\nabla u|^2\phi + k_2(\alpha)\int\alpha u^{\alpha-1}(\Delta u)\phi\leq}\\
\displaystyle{k_1(\alpha)\left[\int\frac{u^{q+\alpha}}{|x|^{\theta}}\phi\right]^{\frac{1}{\chi}}\left[\int|x|^{\frac{\theta\chi'}{\chi}}\frac{|\nabla\phi|^{4\chi'}}{\phi^{3\chi'+\frac{\chi'}{\chi}}}\right]^{\frac{1}{\chi'}}+ \frac{1}{2\epsilon_3^2\alpha}\left[\int\frac{u^{q+\alpha}}{|x|^{\theta}}\phi\right]^{\frac{1}{\chi}}\left[\int|x|^{\frac{\theta\chi'}{\chi}}\frac{|\Delta\phi|^{2\chi'}}{\phi^{\chi'+\frac{\chi'}{\chi}}}\right]^{\frac{1}{\chi'}}} \\+ 
\displaystyle{\lambda\int u^{1+\alpha}V(x)\phi}.\\
\end{array}
\end{equation}
Since the last inequality has only nonnegative terms in the l.h.s., we may let $R\rightarrow+\infty$ and apply (i) together with the Monotone Convergence Theorem in order to proof that $\Delta u\Delta(u^\alpha)\in L^1(\mathbf{R}^N)$ and\\
$$\displaystyle{\int u^{q+\alpha}a(x)+ k_2(\alpha)\int(\Delta u^{\alpha})\Delta u\leq\lambda\int u^{\alpha+1}V(x)}.$$
This last inequality holds true for any small enough positive value of parameters $\epsilon_2$ and $\epsilon_3$. Thus, recalling that $k_2(\alpha)=[1 - \frac{\epsilon_2^2}{2} - \frac{\epsilon_3^2}{2}]$, we let  $\epsilon_2,\epsilon_3\rightarrow 0$ in order to obtain,
\begin{equation*}
\begin{array}{l}
\displaystyle{\int u^{q+\alpha}a(x) + \int(\Delta u^{\alpha})\Delta u \leq \lambda\int u^{\alpha+1}V(x),}
\end{array}
\end{equation*}
as wanted. $\Diamond$
\end{lemma1}
One last tool is required. Its proof can be found in  Theorem 3.9 of \cite{maly}.
\begin{lemma1}\label{lem3.2} Let $r>0$, $u\in W^{2,2}\left(B_{2r}\right)$ be such that $u\geq 0, \Delta u\geq 0$ a.e. in $B_{2r}$. Then 
\begin{equation}\label{2.13} 
||u||_{L^\infty\left({B_r}\right)} \leq C\left(\Mint_{B_{2r}}u^2 dx\right)^\frac{1}{2},\end{equation}
where the constant $C$ does not depend on $r$.\\
\end{lemma1}
Now we can prove the main result.
\begin{theorem1}\label{teo3.1}
Let $u$ be a weak subharmonic solution of (\ref{eq1.1}) in the sense of Definition (\ref{def3.1}), then $u\equiv 0$ a.e.\\
\textbf{Proof.} Choose $\alpha>2$ such that $q<\frac{N-\theta+(4-\theta)\alpha}{N-4}$ and $q+\alpha >4$. If $\gamma = \frac{q+\alpha}{2}$, it is easy to check that $u^\gamma$ satisfies the hypothesis of Lemma \ref{lem3.2} (see Appendix 5.1 for details) $\forall r>0$. Therefore, we have that,
$$||u^\gamma||_{L^\infty\left({B\left(r\right)}\right)}\leq C\left[\Mint_{B(2r)}u^{2\gamma}\right]^\frac{1}{2}\qquad \forall \; r>0,$$
where $C$ is a constant which does not depend on $r$. Moreover, from this last inequality it follows that,
$$\displaystyle{||u^\gamma||_{L^\infty\left({B\left(r\right)}\right)}\leq C\left[\frac{1}{2^Nr^N}\int_{B\left(2r\right)}u^{2\gamma}\right]^\frac{1}{2}=C\left[\frac{1}{r^{N-\theta}}\int_{B\left(2r\right)}\frac{u^{2\gamma}}{r^\theta}\right]^\frac{1}{2}\leq C\left[\frac{1}{r^{N-\theta}}\int_{B\left(2r\right)}\frac{u^{2\gamma}}{|x|^\theta}\right]^\frac{1}{2}},$$
where $C$ does not depend on $r$.\\
Finally, from Lemma \ref{lem3.1} we know that $u^{q+\alpha}a(x)\in L^1(\mathbf{R}^N)$. Therefore, (\ref{eq1.2}) implies  that $\frac{u^{2\gamma}}{|x|^\theta}\in L^1(\mathbf{R}^n)$, and since $N>4>\theta$, we obtain the claim letting $r\rightarrow \infty$. $\; \Diamond$
\end{theorem1}


\section{Appendix}

\subsection{Local Regularity}

In what follows, we review some basic properties of functions in $W^{2,2}_{loc}(\mathbf{R}^N)\bigcap L^\infty(\mathbf{R}^N)$.
For the rest of this section, $\Omega$ will be a regular bounded open set of $\mathbf{R}^N$.\\
First of all, we have the following Gagliardo-Nirenberg interpolation inequality.
\begin{theorem1}\label{teo5.1} Let $u\in L^p(\Omega) \bigcap W^{2,r}(\Omega)$ with $1\leq p\leq\infty$ and $1\leq r\leq\infty$. Then $u\in W^{1,q}(\Omega)$, where $\frac{1}{q}=\frac{1}{2}\left(\frac{1}{p} + \frac{1}{r}\right)$, and
$$||\nabla u||_{L^q} \leq C||u||^\frac{1}{2}_{W^{2,r}}||u||^\frac{1}{2}_{L^p}.$$
\textbf{Proof.} See for instance \cite{nie1}. $\Diamond$
\end{theorem1}
Thanks to this Theorem, it is not hard to show that $W^{2,2}(\Omega)\bigcap L^\infty(\Omega)$ is an algebra. In fact, for any $U,V\in W^{1,2}(\Omega)$, there exist sequences of functions $U_n,V_n\in C^\infty(\mathbf{R}^N)\bigcap W^{2,2}(\Omega)$ such that $U_n\rightarrow U$ and $V_n\rightarrow V$ in $W^{1,2}(\Omega)$ (Theorem 2.3.2 of \cite{ziem}), and this allows us to show that,
 $$\frac{\partial (UV)}{x_i} = V\frac{\partial U}{x_i} + U\frac{\partial V}{x_i} \qquad for \qquad 1\leq i\leq N.$$
Therefore, if $u,v \in W^{2,2}(\Omega)\bigcap L^\infty(\Omega)$, we may apply twice this last equality in order to get that,
\begin{equation}\label{eq4.1} \frac{\partial (uv)}{x_i} = v\frac{\partial u}{x_i} + u\frac{\partial v}{x_i} \;\;\; and \;\;\; \frac{\partial}{\partial x_i}\left(u\frac{\partial v}{\partial x_j}\right) = \frac{\partial u}{\partial x_i}\frac{\partial v}{\partial x_j} + u\frac{\partial^2 v}{\partial x_i\partial x_j} \qquad for \qquad 1\leq i,j\leq N.\end{equation}
The first equality in (\ref{eq4.1}) clearly implies that $uv \in W^{1,2}(\Omega)$ (Thus $W^{1,2}(\Omega)\bigcap L^{\infty}(\Omega)$ is an algebra). Also, applying Theorem \ref{teo5.1} with $p=\infty$ and $r=2$, we get that $\frac{\partial}{\partial x_i}\in L^4(\Omega)$ for $1\leq i\leq N$, and from an applications of Cauchy Schwartz inequality, it follows that,
$$\int_{\Omega}\left|\frac{\partial u}{\partial x_i}\frac{\partial v}{\partial x_j}\right|^2\leq\left[\int_{\Omega}\left|\frac{\partial u}{\partial x_i}\right|^4\right]^2\left[\int_{\Omega}\left|\frac{\partial v}{\partial x_j}\right|^4\right]^2.$$
Therefore, we get from the second equality in (\ref{eq4.1}), that $uv\in W^{2,2}(\Omega)\bigcap L^\infty(\Omega)$ as wanted.\\\\
We also have the following result concerning the composition of a function in a Sobolev space with a function in $C^1(\mathbf{R})$.
\begin{theorem1}\label{teo5.2} Let $G\in C^1(\mathbf{R})$ be such that $G(0)=0\;$ and $|G'(s)|\leq M\;$ $\forall s\in \mathbf{R}$ for some $M>0$. If $u\in W^{1,p}(\Omega)$, then, $G\circ u\in W^{1,p}(\Omega) \;$ and $\; \frac{\partial}{\partial x_i}(G\circ u)= (G'\circ u)\frac{\partial u}{\partial x_i}.$
\textbf{Proof.} See for example Proposition 9.5 of \cite{bz}. $\Diamond$
\end{theorem1}
This theorem can be used to prove that if $u \in W^{2,2}(\Omega)\bigcap L^{\infty}(\Omega)$ and $f\in C^2(\mathbf{R})$, then $f(u)\in W^{2,2}(\Omega)\bigcap L^\infty(\Omega)$ and,
\begin{equation}\label{eq5.2} \frac{\partial^2}{\partial x_i\partial x_j}f(u) = f''(u)\frac{\partial u}{\partial x_i}\frac{\partial u}{\partial x_j} + f'(u)\frac{\partial^2 u}{\partial x_i\partial x_j}. \qquad for \qquad 1\leq i,j\leq N.\end{equation}
In fact, letting
\begin{equation*}G(x) = \begin{cases}f(x)-f(0) & |x| \leq 2||u||_{L^\infty(\Omega)}\\
0 &|x|>3||u||_{L^\infty(\Omega)},\end{cases}\end{equation*}
such that $G\in C^2(\mathbf{R}^N)$, we apply Theorem \ref{teo5.2} in order to get,
$$\frac{\partial}{\partial x_i}\left(\frac{\partial f(u)}{\partial x_j}\right)= \frac{\partial}{\partial x_i}\left(f'(u)\frac{\partial u}{\partial x_j}\right) \qquad for \qquad 1\leq i\leq N,$$
and $f'(u)\in W^{1,2}(\Omega).$ Then, using the fact that $\frac{\partial u}{\partial x_i},f'(u)\in W^{1,2}(\Omega)\;\; for \;\; 1\leq i\leq N$, we obtain (\ref{eq5.2}).\\
The result then follows from the fact that $f'(u), f''(u) \in L^{\infty}(\Omega).$

\end{document}